\documentclass[twoside]{article}
\usepackage{graphicx,amssymb,mathrsfs,amsmath}
\usepackage[numbers,sort&compress]{natbib}
\catcode`@=11
\long\def\@makefntext#1{\noindent #1}
\newskip\tabcentering \tabcentering=1000pt plus 1000pt minus 1000pt
\def\MCH#1#2{\setbox0=\hbox{\raise#1\hbox{#2}}\smash{\box0}}

\def\CR{\cr\noalign{\vspace{1mm} \hrule \vspace{1mm}}}
\def\@evenfoot{}\def\@oddfoot{}

\def\@evenhead{\hbox to\textwidth{\footnotesize\rm\thepage \hfill
{\it Tiexin Guo, Shien Zhao, Xiaolin Zeng}}} 

\def\@oddhead{\hbox to \textwidth{\footnotesize{\it
Random convex analysis (II): continuity and subdifferentiability theorems } \hfill\thepage}}

\def\sec#1{\vskip 3mm\leftline{\bf #1}\vskip 1mm}
\def\subsec#1{\vskip 2mm\leftline{#1}\vskip 1mm}
\def\th#1{\vskip 1mm\noindent{\bf #1}\quad}

\renewcommand{\topfraction}{1}
\renewcommand{\bottomfraction}{1}
\renewcommand{\textfraction}{0}
\renewcommand{\floatpagefraction}{0}
\catcode`@=12

\def\R{{\Bbb R}}

\def\bc{\begin{center}}
\def\ec{\end{center}}
\def\no{\noindent}
\def\hang{\hangindent\parindent}
\def\textindent#1{\indent\llap{\qquad #1\ \ \enspace}\ignorespaces}
\def\ref{\par\hang\textindent}

\def\dl{\displaystyle\lim}

\begin{document}
\abovedisplayskip=6pt plus 1pt minus 1pt \belowdisplayskip=6pt
plus 1pt minus 1pt
\thispagestyle{empty} \vspace*{-1.0truecm} \noindent
\vskip 10mm

\bc{\Large\bf Random convex analysis (II): continuity and subdifferentiability theorems in $L^{0}$--pre--barreled random locally convex modules

\footnotetext{E-mail address:\footnotesize $^{1}$
txguo@buaa.edu.cn,$^{2}$zsefh@cnu.edu.cn, $^{3}$xlinzeng@163.com
}
\footnotetext{\footnotesize Supported by National Natural Science Foundation of China Nos. 11171015, 11301568 and 11401399}} \ec
\bc{\bf Tiexin Guo$^{1}$}\\
{\small School of Mathematics and Statistics, Central South University \\
Changsha {\rm 410083}, P.R. China.}\ec
\bc{\bf Shien Zhao$^{2}$}\\
{\small Elementary Educational College, Capital Normal University\\
Beijing {\rm 100048}, P.R. China.}\ec
\bc{\bf Xiaolin Zeng$^{3}$}\\
{\small School of Mathematics and Statistics, Chongqing Technology and Business University\\
Chongqing {\rm400067}, P.R. China.}\ec

\def\sec#1{\vspace{2mm}\noindent{{\bf #1}}\vspace{0.5mm}}
\def\subsec#1{\vspace{2mm}\leftline{\bf #1}} 
\def\th#1{\vspace{1mm}\noindent{\bf #1}\quad } 
\def\pf#1{\vspace{1mm}\noindent{\it #1}\quad}
\renewcommand{\topfraction}{1} \renewcommand{\bottomfraction}{1}
\renewcommand{\textfraction}{0} \renewcommand{\floatpagefraction}{0}
\floatsep=0pt \textfloatsep=0pt \intextsep=0pt \catcode`@=12
\def\leq{\leqslant}
\def\geq{\geqslant}
\def\R{{\Bbb R}}  \def\N{{\Bbb N}}  \def\Q{{\Bbb Q}}  \def\O{{\Bbb O}} \def\Z{{\Bbb Z}}
\def\C{{\Bbb C}}  \def\hml{\end{document}}  \newsymbol\wjzhml 203F \def\no{\noindent}
\def\CR{\cr\noalign{\vspace{1mm} \hrule \vspace{1mm}}}

\abovedisplayskip=3pt plus 1pt minus 1pt 
\belowdisplayskip=3pt plus 1pt minus 1pt 

\def\le{\leqslant}
\def\ge{\geqslant}
\def\dl{\displaystyle}

\newtheorem{theorem}{Theorem}[section]%
\newtheorem{corollary}[theorem]{Corollary}%
\newtheorem{lemma}[theorem]{Lemma}%
\newtheorem{proposition}[theorem]{Proposition}%
\newtheorem{problem}[theorem]{Problem}
\newtheorem{definition}[theorem]{Definition}%
\newtheorem{remark}[theorem]{Remark}%
\newtheorem{example}[theorem]{Example}




\vspace{8true mm}

\renewcommand{\baselinestretch}{1.9}\baselineskip 19pt

\baselineskip 12pt \renewcommand{\baselinestretch}{1.18}
\noindent{{\bf Abstract}\small\hspace{2.8mm} 
In this paper, we continue to study random convex analysis. First, we introduce the notion of an $L^0$--pre--barreled module. Then, we develop the theory of random duality under the framework of a random locally convex module endowed with the locally $L^0$--convex topology in order to establish a characterization for a random locally convex module to be $L^0$--pre--barreled, in particular we prove that the model space $L^{p}_{\mathcal{F}}(\mathcal{E})$ employed in the module approach to conditional risk measures is $L^0$--pre--barreled, which forms the most difficult part of this paper. Finally, we prove the continuity and subdifferentiability theorems for a proper lower semicontinuous $L^0$--convex function on an $L^{0}$--pre--barreled random locally convex module. So the principal results of this paper may be well suited to the study of continuity and subdifferentiability for $L^0$--convex conditional risk measures.
 }

\vspace{1mm} \no{\footnotesize{\bf Keywords:\hspace{2mm}Random locally convex module,
Locally $L^{0}$--convex topology,
Random conjugate spaces, Proper $L^{0}$--convex lower semicontinuous function, Continuity, Subdifferentiability

}}

\no{\footnotesize{\bf MSC(2000):\hspace{2mm}46A08, 46A20, 46H25,46T20}
 \vspace{2mm}
\baselineskip 15pt
\renewcommand{\baselinestretch}{1.22}
\parindent=10.8pt  
\rm\normalsize\rm

\section{Introduction}

To provide a solid analytic foundation for the module approach to conditional risk measures, in \cite{GZZ} we started to establish a complete random convex analysis, the principal results of \cite{GZZ} are concerned with the study of separation and Fenchel-Moreau duality in random locally convex modules. Based on \cite{GZZ}, this paper continues to study random convex analysis. The main purpose of this paper is to prove continuity and subdifferentiability theorems in $L^{0}$--pre--barreled random locally convex modules.

Continuity and subdifferentiability theorems in classical convex analysis say that a proper lower semicontinuous extended real-valued convex function on a barreled space is continuous and subdifferentiable in the interior of the effective domain of the function. In \cite{FKV}, D. Filipovi\'{c}, M. Kupper and N. Vogelpoth presented the notion of $L^{0}$--barreled modules and established the continuity and subdifferentiability theorems for proper lower semicontinuous $L^{0}$--convex functions defined on $L^{0}$--barreled modules, namely Proposition 3.5 and Theorem 3.7 of \cite{FKV}. Then a natural and key problem is how to characterize an $L^{0}$--barreled module, in particular the problem of whether the model space $L^{p}_{\mathcal{F}}(\mathcal{E})$ employed in the module approach is an $L^{0}$--barreled module or not remains unsolved. Up to now, not even a result characterizing an $L^{0}$--barreled module has been obtained mainly because the notion of an $L^{0}$--barreled module is too similar to that of the classical barreled space and hence also too strong. In this paper, on the basis of our work on separation in \cite{GZZ} we can  overcome the difficulty by presenting the notion of an $L^{0}$--pre-barreled module. The notion of an $L^{0}$--pre-barreled module is weaker than that of an $L^{0}$--barreled module and meets the needs of financial applications. To prove this, we establish random duality theory of a random duality pair under random locally convex modules endowed with the locally $L^{0}$--convex topology so that we can give a characterization for random locally convex modules to be $L^{0}$--pre-barreled, in particular $L^{p}_{\mathcal{F}}(\mathcal{E})$ is $L^{0}$--pre-barreled when it is endowed with the locally $L^{0}$--convex topology, which also forms the most difficult part of this paper. Further, we also prove the new continuity and subdifferentiability theorems based on the notion of an $L^{0}$--pre-barreled module. The results of this paper have been used in \cite{GZZ3}.

In fact, this paper is the second part of our manuscript \cite{GZZ1}. For a random locally convex module $(E,\mathcal{P})$, $\mathcal{P}$ can induce two kinds of topologies, namely the $(\varepsilon,\lambda)$--topology and the locally $L^0$--convex topology. The $(\varepsilon,\lambda)$--topology is very natural in the study of some problems, for example, in \cite{GZZ} we always first consider the related problems under the $(\varepsilon,\lambda)$--topology and then pass to the locally $L^0$--convex topology. On the other hand, the locally $L^0$--convex topology is stronger than the $(\varepsilon,\lambda)$--topology and  often makes some important $L^0$--convex sets possess non-empty interiors, hence this paper employs the locally $L^0$--convex topology for a random locally convex module to meet the needs of continuity and subdifferentiability theorems.

The remainder of this paper is organized as follows. Section 2 is devoted to the study of boundedness of sets in random locally convex modules under locally $L^{0}$--convex topology. Section 3 is devoted to establishing random duality theory of a random duality pair under random locally convex modules endowed with the locally $L^{0}$--convex topology and further characterizing an $L^{0}$--pre-barreled random locally convex module. Finally, in Section 4 we prove the new continuity and subdifferentiability theorems based on the notion of an $L^{0}$--pre-barreled module.

Throughout this paper, we always use the following notation and terminology:

$K:$ the scalar field R of real numbers or C of complex numbers.

$(\Omega,\mathcal{F},P):$ a probability space.

$L^{0}(\mathcal{F},K)=$ the algebra of equivalence classes of $K$--valued $\mathcal{F}$-- measurable random variables
on $(\Omega,\mathcal{F},P)$.

$L^{0}(\mathcal{F})=L^{0}(\mathcal{F},R)$.

$\bar{L}^{0}(\mathcal{F})=$ the set of equivalence classes of extended real-valued $\mathcal{F}$-- measurable random
variables on $(\Omega,\mathcal{F},P)$.

As usual, $\bar{L}^{0}(\mathcal{F})$ is partially ordered by $\xi\leq\eta$ iff $\xi^{0}(\omega)\leq\eta^{0}(\omega)$ for $P$--almost all $\omega\in \Omega$ (briefly, a.s.), where $\xi^0$ and $\eta^0$ are arbitrarily chosen representatives of $\xi$ and $\eta$, respectively. Then $(\bar{L}^{0}(\mathcal{F}),\leq)$ is a complete lattice, $\bigvee H$ and $\bigwedge H$ denote the supremum and infimum of a subset $H$, respectively. $(L^{0}(\mathcal{F}),\leq)$ is a conditionally complete lattice. Please refer to \cite{D-Sch} or \cite[p. 3026]{TXG-JFA} for the rich properties of the supremum and infimum of a set in $\bar{L}^{0}(\mathcal{F})$.

Let $\xi$ and $\eta$ be in $\bar{L}^{0}(\mathcal{F})$. $\xi<\eta$ is understood as usual, namely $\xi\leq\eta$ and $\xi\neq\eta$. In this paper we also use $``\xi<\eta $ (or $\xi\leq\eta$) on $A"$ for $``\xi^{0}(\omega)<\eta^{0}(\omega)$ (resp., $\xi^{0}(\omega)\leq\eta^{0}(\omega)$) for $P$--almost all $\omega\in A"$, where $A\in\mathcal{F}$, $\xi^0$ and $\eta^0$ are arbitrarily chosen representatives of $\xi$ and $\eta$, respectively.

$\bar{L}^{0}_{+}(\mathcal{F})=\{\xi\in \bar{L}^{0}(\mathcal{F})~|~\xi\geq0\}$

$L^{0}_{+}(\mathcal{F})=\{\xi\in L^{0}(\mathcal{F})~|~\xi\geq0\}$

$\bar{L}^{0}_{++}(\mathcal{F})=\{\xi\in \bar{L}^{0}(\mathcal{F})~|~\xi>0$ on $\Omega\}$

$L^{0}_{++}(\mathcal{F})=\{\xi\in L^{0}(\mathcal{F})~|~\xi>0$ on $\Omega\}$

Besides, $\tilde{I}_{A}$ always denotes the equivalence class of $I_{A}$, where $A\in \mathcal{F}$ and $I_{A}$ is the characteristic function of $A$.
When $\tilde{A}$ denotes the equivalence class of $A (\in \mathcal{F})$, namely $\tilde{A}=\{B\in\mathcal{F}~|~P(A\triangle B)=0\}$ (here, $A\triangle B=(A\setminus B)\bigcup(B\setminus A)$), we also use $I_{\tilde{A}}$ for $\tilde{I}_{A}$.

Specially, $[\xi<\eta]$ denotes the equivalence class of $\{\omega\in\Omega~|~\xi^0(\omega)<\eta^0(\omega)\}$, where $\xi^0$ and $\eta^0$ are arbitrarily chosen representatives of $\xi$ and $\eta$ in $\bar{L}^0(\mathcal{F})$, respectively, some more notations such as $[\xi=\eta]$ and $[\xi\neq\eta]$ can be similarly understood.

\section{Boundedness of sets in random locally convex modules under locally $L^{0}$--convex topology}

The main results in this section are Theorems 2.13 and 2.14 below. Let us first recall some basic notions and terminology.

\begin{definition}($See$ \cite{TXG-Master,TXG-PHD,TXG-basic}). An ordered pair $(E,\|\cdot\|)$ is called a random normed space (briefly, an $RN$ space) over $K$ with base $(\Omega,\mathcal{F},P)$ if $E$ is a linear space over $K$ and $\|\cdot\|$ is a mapping from $E$ to $L^{0}_{+}(\mathcal{F})$ such that the following are satisfied:

\noindent ($RN$--1). $\|\alpha x\|=|\alpha| \|x\|$, $\forall \alpha\in K$ and $x\in E$;

\noindent ($RN$--2). $\|x\|=0$ implies $x=\theta$ (the null element of $E$);

\noindent ($RN$--3). $\|x+y\|\leq\|x\|+\|y\|$, $\forall x,y\in E$.

\noindent Here $\|\cdot\|$ is called the random norm on $E$ and $\|x\|$ the random norm of $x\in E$ (If $\|\cdot\|$ only satisfies ($RN$--1) and ($RN$--3) above, it is called a random seminorm on $E$).

\noindent Furthermore, if, in addition, $E$ is a left module over the algebra $L^{0}(\mathcal{F},K)$ (briefly, an $L^{0}(\mathcal{F},K)$--module) such that

\noindent ($RNM$--1). $\|\xi x\|=|\xi| \|x\|$, $\forall \xi\in L^{0}(\mathcal{F},K)$ and $x\in E$.

\noindent Then $(E,\|\cdot\|)$ is called a random normed module (briefly, an $RN$ module) over $K$ with base $(\Omega,\mathcal{F},P)$, the random norm $\|\cdot\|$ with the property ($RNM$--1) is also called an $L^0$--norm on $E$ (a mapping only satisfying ($RN$--3) and ($RNM$--1) above is called an $L^0$--seminorm on E).

\end{definition}

\begin{definition}($See$ \cite{TXG-Master,TXG-PHD,TXG-Module,TXG-basic}). Let $(E_1,\|\cdot\|)$ and $(E_2,\|\cdot\|)$ be $RN$ spaces over $K$ with base $(\Omega,\mathcal{F},P)$. A linear operator $T$ from $E_1$ to $E_2$ is said to be a.s. bounded if there is $\xi\in L^{0}_{+}(\mathcal{F})$ such that $\|Tx\|_2\leq\xi\|x\|_1, \forall x\in E_1$. Denote by $B(E_1,E_2)$ the linear space of a.s. bounded linear operators from $E_1$ to $E_2$, define $\|\cdot\|:B(E_1,E_2)\rightarrow L^{0}_{+}(\mathcal{F})$ by $\|T\|=\bigwedge\{\xi\in L^{0}_{+}(\mathcal{F})~|~\|Tx\|_2\leq\xi\|x\|_1$ for all $x\in E_1\}$ for all $T\in B(E_1,E_2)$, then it is easy to check that $(B(E_1,E_2),\|\cdot\|)$ is also an $RN$ space over $K$ with base $(\Omega,\mathcal{F},P)$, in particular $(B(E_1,E_2),\|\cdot\|)$ is an $RN$ module if so is $E_2$. Specially, for a fixed random normed space $(E,\|\cdot\|)$ over $K$ with base $(\Omega,\mathcal{F},P)$, the $RN$ module $(E^{\ast},\|\cdot\|)$ with $E^{\ast}=B(E,L^{0}(\mathcal{F},K))$ is called the random conjugate space of $E$, where $L^{0}(\mathcal{F},K)$ is always regarded as an $RN$ module endowed with the $L^0$--norm $|\cdot|$.
\end{definition}

\begin{definition}($See$ \cite{TXG-PHD,TXG-Module,TXG-Sur}). An ordered pair $(E,\mathcal{P})$ is called a random locally convex space (briefly, an $RLC$ space) over $K$ with base $(\Omega,\mathcal{F},P)$ if $E$ is a linear space over $K$ and $\mathcal{P}$ a family of mappings from $E$ to $L^0_{+}(\mathcal{F})$ such that the following are satisfied:

\noindent ($RLC$--1). Every $\|\cdot\|\in \mathcal{P}$ is a random seminorm on $E$;

\noindent ($RLC$--2). $\bigvee\{\|x\|:\|\cdot\|\in\mathcal{P}\}=0$ iff $x=\theta$.

\noindent Furthermore, if, in addition, $E$ is an $L^{0}(\mathcal{F},K)$--module and each $\|\cdot\|\in \mathcal{P}$ is an $L^0$--seminorm on $E$, then $(E,\mathcal{P})$ is called a random locally convex module (briefly, an $RLC$ module) over $K$ with base $(\Omega,\mathcal{F},P)$.
\end{definition}

In the sequel of this paper, given a random locally convex space $(E,\mathcal{P})$, $\mathcal{P}_f$ always denotes the family of finite subsets of $\mathcal{P}$ and for each $\mathcal{Q}\in \mathcal{P}_f$, $\|\cdot\|_{\mathcal{Q}}$ denotes the random seminorm defined by $\|x\|_{\mathcal{Q}}=\bigvee\{\|x\|:\|\cdot\|\in \mathcal{Q}\}$ for all $x\in E$.

Following is an important example used in this paper.

D. Filipovi$\acute{c}$, M. Kupper and N. Vogelpoth constructed important $RN$ modules $L^p_{\mathcal{F}}(\mathcal{E}) (1\leq p\leq+\infty)$ in \cite{FKV}.

\begin{example} Let $(\Omega, {\mathcal E}, P)$ be a probability space and ${\mathcal F}$
a sub--$\sigma$--algebra of ${\mathcal E}$. Define $|||\cdot|||_p\colon L^0({\mathcal E})\to {\bar L}^0_+({\mathcal F})$ by
$$|||x|||_p=\left\{
               \begin{array}{ll}
                 E[|x|^p|{\mathcal F}]^{1\over p}, & \hbox{when $1\leq p<\infty$;} \\
                 \bigwedge\{\xi\in {\bar L}^0_+({\mathcal F})~|~|x|\leq\xi\}, & \hbox{when $p=+\infty$;}
               \end{array}
             \right.
$$
for all $x\in L^0({\mathcal E})$.

Denote $L^p_{\mathcal F}({\mathcal E})=\{x\in L^0({\mathcal E})~|~|||x|||_p\in L^0_+({\mathcal F})\}$, then $(L^p_{\mathcal F}({\mathcal E}), |||\cdot|||_p)$ is an $RN$ module over $R$ with base $(\Omega, {\mathcal F}, P)$ and $L^p_{\mathcal F}({\mathcal E})=L^0({\mathcal F})\cdot L^p({\mathcal E})=\{~\xi x~|~\xi\in L^0({\mathcal F})~ \hbox{and}~ x\in L^p({\mathcal E})\}$.

\end{example}

\begin{definition}($See$ \cite{TXG-PHD,TXG-Module,TXG-Sur,TXG-SLP}). Let $(E, {\mathcal P})$ be an $RLC$ space over $K$ with base $(\Omega, {\mathcal F}, P)$. For any positive numbers $\varepsilon$ and $\lambda$ with $0<\lambda<1$ and $\mathcal{Q}\in {\mathcal P}_f$, let $N_{\theta}(\mathcal{Q}, \varepsilon, \lambda)=\{x\in E~|~P\{\omega\in \Omega~|~\|x\|_\mathcal{Q}(\omega)<\varepsilon\}>1-\lambda\}$, then $\{N_{\theta}(\mathcal{Q}, \varepsilon, \lambda)~|~\mathcal{Q}\in {\mathcal P}_f, \varepsilon >0, 0<\lambda<1\}$ forms a local base at $\theta$ of some Hausdorff linear topology on $E$, called the $(\varepsilon, \lambda)$--topology induced by ${\mathcal P}$.
\end{definition}

From now on, we always denote by ${\mathcal T}_{\varepsilon, \lambda}$ the $(\varepsilon, \lambda)$--topology for every $RLC$ space if there is no possible confusion. Clearly, the $(\varepsilon, \lambda)$--topology for the special $RN$ module $L^0({\mathcal F}, K)$ is exactly the ordinary topology of convergence in measure, and $(L^0({\mathcal F}, K), {\mathcal T}_{\varepsilon, \lambda})$ is a topological algebra over $K$. It is also easy to check that $(E, {\mathcal T}_{\varepsilon, \lambda})$ is a topological module over $(L^0({\mathcal F}, K), {\mathcal T}_{\varepsilon, \lambda})$ when $(E, {\mathcal P})$ is an $RLC$ module over $K$ with base $(\Omega, {\mathcal F}, P)$, namely the module multiplication operation is jointly continuous.

For an $RLC$ module $(E, {\mathcal P})$ over $K$ with base $(\Omega, {\mathcal F}, P)$, we always denote by $(E, {\mathcal P})^\ast_{\varepsilon, \lambda}$ ( or, briefly, $E^\ast_{\varepsilon, \lambda}$, whenever there is no confusion ) the $L^0({\mathcal F}, K)$--module of continuous module homomorphisms from $(E, {\mathcal T}_{\varepsilon, \lambda})$ to $(L^0({\mathcal F}, K), {\mathcal T}_{\varepsilon, \lambda})$, called the random conjugate space of $(E, {\mathcal P})$ under the  $(\varepsilon, \lambda)$--topology.

\begin{proposition}($See$ \cite{TXG-PHD,TXG-Module,Guotx-onsome}). Let $(E_1, \|\cdot\|_1)$ and $(E_2, \|\cdot\|_2)$ be two $RN$ modules over $K$ with base $(\Omega, {\mathcal F}, P)$ and $T$ a linear operator from $E_1$ to $E_2$. Then $T\in B(E_1, E_2)$ iff $T$ is a continuous module homomorphism from $(E_1, {\mathcal T}_{\varepsilon, \lambda})$ to $(E_2, {\mathcal T}_{\varepsilon, \lambda})$, in which case $\|T\|=\bigvee\{\|Tx\|_2~|~x\in E_1~\hbox{and}~ \|x\|_1\leq 1\}$.
\end{proposition}

Proposition 2.6 is very useful, Guo uses it to prove that $(B(E_1, E_2), \|\cdot\|)$ is always ${\mathcal T}_{\varepsilon, \lambda}$--complete for any two $RN$ spaces $E_1$ and $E_2$ such that $E_2$ is ${\mathcal T}_{\varepsilon, \lambda}$--complete, in particular $E^\ast$ is ${\mathcal T}_{\varepsilon, \lambda}$--complete for every $RN$ space $E$, cf. \cite{TXG-PHD,TXG-Module}. It is also clear from Proposition 2.6 that $E^\ast=E^\ast_{\varepsilon, \lambda}$ for every $RN$ module $E$, cf. \cite{TXG-PHD,TXG-Extension}.

For any $\varepsilon \in L^0_{++}({\mathcal F})$, let $U(\varepsilon)=\{\xi\in L^0({\mathcal F}, K)~|~|\xi|\leq \varepsilon\}$. A subset $G$ of $L^0({\mathcal F}, K)$ is ${\mathcal T}_c$--open if for each fixed $x\in G$ there is some $\varepsilon \in L^0_{++}({\mathcal F})$ such that $x+U(\varepsilon)\subset G$. Denote by ${\mathcal T}_c$ the family of ${\mathcal T}_c$--open subsets of $L^0({\mathcal F}, K)$, then ${\mathcal T}_c$ is a Hausdorff topology on $L^0({\mathcal F}, K)$ such that $(L^0({\mathcal F}, K), {\mathcal T}_c)$ is a topological ring, namely the addition and multiplication operations are jointly continuous. D. Filipovi\'{c}, M. Kupper and N. Vogelpoth first observed this kind of topology and further pointed out that ${\mathcal T}_c$ is not necessarily a linear topology since the mapping $\alpha\mapsto \alpha x$ ($x$ is fixed) is no longer continuous in general. These observations led them to the study of a class of topological modules over the topological ring $(L^0({\mathcal F}, K), {\mathcal T}_c)$ in \cite{FKV}, where they only considered the case when $K=R$, in fact the complex case can also similarly introduced as follows.

\begin{definition}($See$ \cite{FKV}). An ordered pair $(E, {\mathcal T})$ is a topological $L^0({\mathcal F}, K)$--module if both $(E, {\mathcal T})$ is a topological space and $E$ is an $L^0({\mathcal F}, K)$--module such that $(E, {\mathcal T})$ is a topological module over the topological ring $(L^0({\mathcal F}, K), {\mathcal T}_c)$, namely the addition and module multiplication operations are jointly continuous.
\end{definition}

Denote by $(E, {\mathcal T})^\ast$ ( briefly, $E^\ast_c$ ) the $L^0({\mathcal F}, K)$--module of continuous module homomorphisms from $(E, {\mathcal T})$ to $(L^0({\mathcal F}, K), {\mathcal T}_c)$, called the random conjugate space of the topological $L^0({\mathcal F}, K)$--module $(E, {\mathcal T})$, which was first introduced in \cite{FKV}.

\begin{definition}($See$ \cite{TXG-Sur,TXG-XXC,FKV}). Let $E$ be an $L^0({\mathcal F}, K)$--module and $A$ and $B$ two subsets of $E$. $A$ is said to be $L^0$--absorbed by $B$ if there is some $\xi \in L^0_{++}({\mathcal F})$ such that $\eta A\subset B$ for all $\eta\in L^0({\mathcal F}, K)$ with $|\eta|\leq \xi$. $B$ is $L^0$--absorbent if $B$ $L^0$--absorbs every element in $E$. $B$ is $L^0$--convex if $\xi x+(1-\xi)y\in B$ for all $x,\,y\in B$ and $\xi\in L^0_+({\mathcal F})$ with $0\leq \xi\leq 1$. $B$ is $L^0$--balanced if $\eta B\subset B$ for all $\eta\in L^0({\mathcal F}, K)$ with $|\eta|\leq 1$.
\end{definition}

\begin{definition}($See$ \cite{FKV}). A topological $L^0({\mathcal F}, K)$--module $(E, {\mathcal T})$ is called a locally $L^0$--convex $L^0({\mathcal F}, K)$--module ( briefly, a locally $L^0$--convex module when $K=R$ ), in which case ${\mathcal T}$ is called a locally $L^0$--convex topology on $E$, if ${\mathcal T}$ has a local base ${\mathcal B}$ at $\theta$ ( the null element in $E$ ) such that each member in ${\mathcal B}$ is $L^0$--balanced, $L^0$--absorbent and $L^0$--convex.
\end{definition}

\begin{proposition}($See$ \cite{FKV}). Let ${\mathcal P}$ be a family of $L^0$--seminorms on an $L^0({\mathcal F}, K)$--module $E$. For any $\varepsilon \in L^0_{++}({\mathcal F})$ and any $ Q\in {\mathcal P}_f$ (namely $Q$  is a finite subset of ${\mathcal P}$), let $N_{\theta}(Q, \varepsilon)=\{x\in E~|~\|x\|_Q\leq \varepsilon\}$, then $\{~N_{\theta}(Q, \varepsilon)~|~Q\in {\mathcal P}_f, ~\varepsilon \in L^0_{++}({\mathcal F})\}$ forms a local base at $\theta$ of some locally $L^0$--convex topology, called the locally $L^0$--convex topology induced by ${\mathcal P}$. Specially, Let $(E,{\mathcal P})$ be an $RLC$ module over $K$ with base $(\Omega, {\mathcal F}, P)$ and ${\mathcal T}_c$ the locally $L^0$--convex topology induced by ${\mathcal P}$. Then $(E, {\mathcal T}_c)$ is a Hausdorff locally $L^0$--convex $L^0({\mathcal F}, K)$--module.
\end{proposition}

From now on, we always denote by $\mathcal{T}_c$ the locally $L^0$--convex topology induced by $\mathcal{P}$ for every $RLC$ module $(E,\mathcal{P})$ if there is no risk of confusion.

Let $(E,{\mathcal P})^\ast_c=(E,{\mathcal T}_c)^\ast$ (briefly, $E^\ast_c$, if there is no risk of confusion), called the random conjugate space of a random locally convex module $(E,{\mathcal P})$ under the locally $L^0$--convex topology ${\mathcal T}_c$ induced by ${\mathcal P}$.

Let $(E,\mathcal{P})$ be a random locally convex space. Since the $(\varepsilon,\lambda)$--topology induced by $\mathcal{P}$ is a linear topology, the notion of boundness under the $(\varepsilon,\lambda)$--topology of a set in $E$ is as usual, this kind of bounded sets are important in some fields, for example, they are called ``probabilistically bounded sets" in the theory of probabilistic normed spaces (see \cite{SS}) and are called ``stochastically bounded sets" for Banach space-valued random elements in probability theory in Banach spaces, whereas the following notion of bounded sets will play a crucial role in random duality theory in this paper as well as in \cite{TXG-XXC}.

\begin{definition}
Let $(E,\mathcal{T})$ be a locally $L^{0}$--convex $L^{0}(\mathcal{F},K)$--module. $A\subset E$ is said to be $\mathcal{T}$--bounded if $A$ can be $L^0$--absorbed by every neighborhood of $\theta$.

\end{definition}

From now on, we always suppose that all the $L^0({\mathcal F}, K)$--modules $E$ involved in this paper have the property that for any $x,~y\in E$, if there is a countable partition $\{A_n,n\in N\}$ of $\Omega$ to ${\mathcal F}$ such that ${\tilde I}_{A_n}x={\tilde I}_{A_n}y$ for each $n\in N$ then $x=y$. Guo already pointed out in \cite{TXG-JFA} that all random locally convex modules possess this property, so the assumption is not too restrictive.

\begin{definition}($See$ \cite{TXG-JFA}). Let $E$ be an $L^0({\mathcal F}, K)$--module. A sequence $\{x_n, n\in N\}$ in $E$ is countably concatenated in $E$ with respect to a countable partition $\{A_n,n\in N\}$ of $\Omega$ to ${\mathcal F}$ if there is $x\in E$ such that ${\tilde I}_{A_n}x={\tilde I}_{A_n}x_n$ for each $n\in N$, in which case we define $\sum^{\infty}_{n=1}{\tilde I}_{A_n}x_n$ as $x$. A subset $G$ of $E$ is said to have the countable concatenation property if each sequence $\{x_n, n\in N\}$ in $G$ is countably concatenated in $E$ with respect to an arbitrary countable partition $\{A_n,n\in N\}$ of $\Omega$ to $\mathcal{F}$ and $\sum^{\infty}_{n=1}{\tilde I}_{A_n}x_n\in G$.
\end{definition}

Let $E$ be an $L^0({\mathcal F}, K)$--module with the countable concatenation property, from now on for a subset $G$ of $E$ we always use $H_{cc}(G)$ for the countable concatenation hull of $G$, namely $H_{cc}(G)=\{\Sigma_{n=1}^{\infty}\tilde{I}_{A_n}x_n:\{x_n,n\in N\}$ is a sequence in $G$ and $\{A_n,n\in N\}$ is a countable partition of $\Omega$ to $\mathcal{F}\}$

We can now state the main results in this section.

\begin{theorem} $(${\bf Resonance theorem}$)$ Let $(E_1,\|\cdot\|_1)$ and $(E_2,\|\cdot\|_2)$ be two $RN$ modules over $K$ with base $(\Omega,\mathcal{F},P)$ such that $E_1$ is $\mathcal{T}_c$--complete and has the countable concatenation property. For a subset $\{T_{\alpha},\alpha\in\Lambda\}$ of $B(E_1,E_2)$, then $\{T_{\alpha},\alpha\in\Lambda\}$ is $\mathcal{T}_c$--bounded in $B(E_1,E_2)$ iff $\{T_{\alpha}x,\alpha\in\Lambda\}$ is $\mathcal{T}_c$--bounded in $E_2$ for all $x\in E_1$.
\end{theorem}

Let $(E_1,\|\cdot\|)$ and $(E_2,\|\cdot\|)$ be $RN$ modules over $K$ with base $(\Omega,\mathcal{F},P)$. It is easy to prove that a linear operator $T:E_1\rightarrow E_2$ belongs to $B(E_1,E_2)$ iff $T$ is a continuous module homomorphism from $(E_1,\mathcal{T}_c)$ to $(E_2,\mathcal{T}_c)$. Hence Theorem 2.13 also gives a resonance theorem for a family of continuous module homomorphisms from $(E_1,\mathcal{T}_c)$ to $(E_2,\mathcal{T}_c)$.

\begin{theorem} Let $(E,\mathcal{P})$ be an $RLC$ module over $K$ with base $(\Omega,\mathcal{F},P)$ and $A\subset E$. Then $A$ is $\mathcal{T}_c$--bounded iff $f(A)$ is $\mathcal{T}_c$--bounded in $(L^{0}(\mathcal{F},K),\mathcal{T}_c)$ for every $f\in E^{\ast}_c$.
\end{theorem}

Theorems 2.13 and 2.14 are implied by the work on resonance theorem at the earlier stage of $RLC$ spaces. To see this, let us recall:

\begin{definition} (See \cite{TXG-PHD,TXG-Module}). Let $(E,\mathcal{P})$ be an $RLC$ space over $K$ with base $(\Omega,\mathcal{F},P)$. A set $A\subset E$ is said to be a.s. bounded if $\bigvee\{\|a\|:a\in A\}\in L^{0}_{+}(\mathcal{F})$ for each $\|\cdot\|\in\mathcal{P}$.
\end{definition}

Lemma 2.16 below is clear by definition.

\begin{lemma} Let $(E,\mathcal{P})$ be an $RLC$ module. Then a set $A$ of $E$ is $\mathcal{T}_c$--bounded iff $A$ is a.s. bounded.
\end{lemma}

The proof of Theorem 2.13 remains to need Propositions 2.17 and 2.18 below.

\begin{proposition}($See$ \cite{TXG-JFA}). Let $(E, {\mathcal P})$ be an $RLC$ module. Then $E$ is ${\mathcal T}_{\varepsilon, \lambda}$--complete iff both $E$ has the countable concatenation property and $E$ is ${\mathcal T}_c$--complete.
\end{proposition}

\begin{proposition} ($See$ \cite{TXG-PHD,TXG-Module,Guotx-onsome}). Let $(E_1,\|\cdot\|_1)$ and $(E_2,\|\cdot\|_2)$ be two $RN$ modules over $K$ with base $(\Omega,\mathcal{F},P)$ such that $E_1$ is $\mathcal{T}_{\varepsilon,\lambda}$--complete. Given a subset $\{T_{\alpha},\alpha\in\Lambda\}$ in $B(E_1,E_2)$, then $\{T_{\alpha},\alpha\in\Lambda\}$ is a.s. bounded in $B(E_1,E_2)$ iff  $\{T_{\alpha}x,\alpha\in\Lambda\}$ is a.s. bounded in $E_2$ for all
$x\in E_1$.
\end{proposition}

We can now prove Theorem 2.13.

\noindent{\bf Proof of Theorem 2.13.}
It immediately follows from Lemma 2.16, Propositions 2.17 and 2.18. \hfill $\square$

For the proof of Theorem 2.14, let us recall from \cite{TXG-PHD,TXG-Module}: Let $(E,\mathcal{P})$ be an $RLC$ space over $K$ with base $(\Omega,\mathcal{F},P)$. A linear operator $f$ from $E$ to $L^0(\mathcal{F},K)$ (such an operator is also called a random linear functional on $E$) is called an a.s. bounded random linear functional of type I if there are $\xi\in L^0_+(\mathcal{F})$ and some finite subset $\mathcal{Q}$ of $\mathcal{P}$ such that $|f(x)|\leq\xi \|x\|_{\mathcal{Q}}$ for all $x\in E$. Denote by $E^{\ast}_{I}$ the $L^0(\mathcal{F},K)$--module of a.s. bounded random linear functionals on $E$ of type I, called the first kind of random conjugate space of $(E,\mathcal{P})$.

The proof of Theorem 2.14 remains to need Propositions 2.19 and 2.20 below.

\begin{proposition} ($See$ \cite{TXG-PHD,TXG-Module,TXG-Sur}). Let $(E,\mathcal{P})$ be an $RLC$ space over $K$ with base $(\Omega,\mathcal{F},P)$ and $A$ a subset of $E$. Then $A$ is a.s. bounded iff $f(A)$ is a.s. bounded in $(L^0(\mathcal{F},K),|\cdot|)$ for each $f\in E^{\ast}_{I}$.

\end{proposition}

\begin{proposition}($See$ \cite{GZZ}). Let $(E,{\mathcal P})$ be a random locally convex module over $K$ with base $(\Omega, {\mathcal F}, P)$ and $f\colon E\to L^0({\mathcal F}, K)$ a random linear functional. Then $f\in E^\ast_I$ iff $f\in E^\ast_c$, namely $E^\ast_I=E^\ast_c$.
\end{proposition}

We can now prove Theorem 2.14.

\noindent{\bf Proof of Theorem 2.14.}
It immediately follows from Propositions 2.19 and 2.20. \hfill $\square$

\section{Random duality under the locally $L^0$--convex topology with respect to random duality pair}

Only the classical duality theory with respect to a duality pair can give a thorough treatment of classical conjugate space theory of locally convex spaces, cf. \cite{DXX}. The theory of random conjugate spaces  occupies a central place in the study of $RN$ modules and $RLC$ modules, it is very natural that random duality theory was studied at the previous time in \cite{TXG-dual,TXG-Sur,TXG-XXC}, where many basic results and useful techniques were already obtained. Before 2009, only the $(\varepsilon, \lambda)$--topology was available, so the work in \cite{TXG-dual,TXG-Sur,TXG-XXC} was carried out under this topology, where the family of $L^0$--seminorms plays a key role. In this section, we will establish some basic results on random duality theory with respect to the locally $L^0$--convex topology in order to provide an enough framework for the theory of $RLC$ modules and its financial applications.

\subsection{Random compatible locally $L^0$--convex topology}

Let $X$ and $Y$ be two $L^0({\mathcal F}, K)$--modules. A mapping $\langle\cdot, \cdot\rangle\colon X\times Y\to L^0({\mathcal F}, K)$ is said to be $L^0({\mathcal F}, K)$--bilinear function ( briefly, $L^0$--bilinear function if there is not any possible confusion ) if both $\langle x, \cdot\rangle\colon Y\to L^0({\mathcal F}, K)$ and $\langle\cdot, y\rangle\colon X\to L^0({\mathcal F}, K)$ are $L^0$--linear functions for all $x\in X$ and $y\in Y$.

\begin{definition} (See \cite{TXG-dual,TXG-Sur,TXG-XXC}). Two $L^0({\mathcal F}, K)$--modules $X$ and $Y$ are called a random duality pair over $K$ with base $(\Omega, {\mathcal F}, P)$ with respect to the $L^0$--bilinear function $\langle\cdot,\cdot\rangle\colon X\times Y\to L^0({\mathcal F}, K)$ if the following are satisfied:\\
(1). $\langle x, y\rangle=0$ for all $y\in Y$ iff $x=\theta$;\\
(2). $\langle x, y\rangle=0$ for all $x\in X$ iff $y=\theta$.
\end{definition}

Usually, if $X,~Y$ and $\langle\cdot,\cdot\rangle$ satisfy Definition 3.19, then we simply say that $\langle X, Y\rangle$ is a random duality pair over $K$ with base $(\Omega, {\mathcal F}, P)$. Let $X^\#$ denote the $L^0({\mathcal F}, K)$--module of $L^0$--linear functions from an $L^0({\mathcal F}, K)$--module $X$ to $L^0({\mathcal F}, K)$. It is clear that $X^\#$ has the countable concatenation property. If $\langle X, Y\rangle$ is a random duality pair, we always identify each $x\in X$ with $\langle x, \cdot\rangle\in Y^\#$, namely regard $X$ as a submodule of $Y^\#$, thus for any subset $G\subset X$, we always use $H_{cc}(G)$ for the countable concatenation hull of $G$ in $Y^\#$, which would not cause any possible confusion.

\begin{definition} Let $\langle X, Y\rangle$ be a random duality pair over $K$ with base $(\Omega, {\mathcal F}, P)$. A family ${\mathcal P}$ of $L^0$--seminorms on $X$ is called a random compatible family with $Y$ with respect to the locally $L^0$--convex topology ${\mathcal T}_c$ induced by ${\mathcal P}$ if $(X, {\mathcal P})$ becomes an $RLC$ module over $K$ with base $(\Omega, {\mathcal F}, P)$ such that $(E, {\mathcal P})^\ast_c=Y$, in which case we also say that ${\mathcal T}_c$ is a random compatible locally $L^0$--convex topology with $Y$.
\end{definition}

\begin{remark} In \cite{TXG-XXC}, a family ${\mathcal P}$ of $L^0$--seminorms is called a random compatible family with $Y$ with respect to the $(\varepsilon, \lambda)$--topology ${\mathcal T}_{\varepsilon, \lambda}$ induced by ${\mathcal P}$ if $(X, {\mathcal P})$ becomes an $RLC$ module such that $(E, {\mathcal P})^\ast_{\varepsilon, \lambda}=Y$.
\end{remark}

Let $\langle X, Y\rangle$ be a random duality pair and $\sigma(X, Y)=\{|\langle\cdot, y\rangle|\colon y\in Y\}$, then $\sigma(X, Y)$ is a family of $L^0$--seminorms on $X$ such that $(X, \sigma(X, Y))$ becomes an $RLC$ module. In the sequel, $\sigma_c(X, Y)$ and $\sigma_{\varepsilon, \lambda}(X, Y)$ always denote the locally $L^0$--convex topology and the $(\varepsilon, \lambda)$--topology induced by $\sigma(X, Y)$, respectively.

\begin{theorem} Let $\langle X, Y\rangle$ be a random duality pair over $K$ with base $(\Omega, {\mathcal F}, P)$. Then $\sigma_c(X, Y)$ is a random compatible topology with $Y$.
\end{theorem}

Proof of Theorem 3.4 needs Lemma 3.5 below.

\begin{lemma} $($See \cite{TXG-dual,TXG-XXC}$.)$ Let $E$ be an $L^0({\mathcal F}, K)$--module. If $f_1,f_2,\cdots,f_n$ and $g$ are $n+1$ $L^0$--linear functions from $E$ to $L^0({\mathcal F}, K)$, then there are $\xi_1,\xi_2,\cdots,\xi_n\in L^0({\mathcal F}, K)$ such that $g=\sum^n_{i=1}\xi_if_i$ iff $\bigcap^n_{i=1}N(f_i)\subset N(g)$, where $N(f_i)=\{x\in E~|~f_i(x)=0\}$ ( $1\leq i\leq n$ ) and $N(g)=\{x\in E~|~g(x)=0\}$.
\end{lemma}

We can now prove Theorem 3.4.

\noindent{\bf Proof of Theorem 3.4.}
Since it is obvious that $Y\subset (X, \sigma(X, Y))^\ast_c$, it remains to prove that $(X, \sigma(X, Y))^\ast_c \subset Y$. Let $f\in (X, \sigma(X, Y))^\ast_c$, then by Proposition 2.20 there are $\xi\in L^0_+({\mathcal F})$ and $y_1,y_2,\cdots, y_n\in Y$ such that $|f(x)|\leq \xi (\bigvee\{|\langle x, y_i\rangle|\colon 1\leq i\leq n\})$ for all $x\in X$. By Lemma 3.5, there are $\xi_1,\xi_2,\cdots,\xi_n\in L^0({\mathcal F}, K)$ such that $f=\sum^n_{i=1}\xi_i\langle x, y_i\rangle$ for all $x\in X$. Let $y=\sum^n_{i=1}\xi_i y_i$, then $f=y\in Y$.  \hfill $\square$

\begin{remark} In \cite{TXG-XXC}, since we employed the $(\varepsilon, \lambda)$--topology, we proved that $(X,\sigma(X, Y))^\ast_{\varepsilon, \lambda}=H_{cc}(Y)$, which motivates us to find out Theorem 3.7 of \cite{GZZ}. In \cite{TXG-XXC}, $Y$ is regular with respect to $X$ if, for each sequence $\{y_n, n\in N\}$ and each countable partition $\{A_n, n\in N\}$ of $\Omega $ to ${\mathcal F}$, there is $y\in Y$ such that $\langle x, y\rangle=\sum^{\infty}_{n=1}{\tilde I}_{A_n}\langle x, y_n\rangle$ for all $x\in X$, which implies that ${\tilde I}_{A_n}\langle x, y\rangle={\tilde I}_{A_n}\langle x, y_n\rangle$ for all $x\in X$ and $n\in N$, namely ${\tilde I}_{A_n}y={\tilde I}_{A_n}y_n$ for each $n\in N$, that is to say, $Y$ has the countable concatenation property. Thus, for a random duality pair $\langle X, Y\rangle$, `` $Y$ is regular '' and `` $Y$ has the countable concatenation property '' are the same thing.
\end{remark}

For the proof of Theorem 3.8, we need Lemma 3.7 below.

\begin{lemma} ($See$ \cite{TXG-XXC}.) Let $X$ be an $L^{0}(\mathcal{F},K)$--module, $f:X\rightarrow L^{0}(\mathcal{F},K)$ an $L^{0}(\mathcal{F},K)$--linear function and $\{p_i:X\rightarrow L^{0}_+(\mathcal{F})~|~1\leq i\leq n\}$ n $L^{0}$--seminorms such that $|f(x)|\leq\Sigma_{i=1}^{n}p_i(x)$ for all $x\in X$. Then for each $1\leq i\leq n$ there is an $L^{0}(\mathcal{F},K)$--linear function $f_i$ such that

\noindent (1). $|f_i(x)|\leq p_i(x)$ for all $x\in X$  and $1\leq i\leq n$;

\noindent (2). $f(x)=\Sigma_{i=1}^{n}f_n(x)$ for all $x\in X$.
\end{lemma}

\begin{theorem} Let $\langle X, Y\rangle$ be a random duality pair. Then there is a strongest one in all the random compatible locally $L^0$--convex topologies with $Y$.
\end{theorem}

\begin{proof} By Lemma 3.7, one can see that the proof is completely similar to the one of the corresponding classical case, so is omitted.
\end{proof}

Let $\langle X, Y\rangle$ be a random duality pair. For a subset $A$ of $X$, $A^0:=\{y\in Y|~|\langle a, y\rangle|\leq 1~\hbox{for all $a\in A$~}\}$ is called the polar of $A$ in $Y$. Similarly, one can define the polar of a subset $B$ of $Y$ in $X$. For the proof of Theorem 3.11 below, we need Lemmas 3.9 and 3.10 below.

\begin{lemma}($See$ \cite{TXG-JFA}.) Let $(E, {\mathcal P})$ be an $RLC$ module and $G$ a subset of $E$ such that $G$ has the countable concatenation property. Then ${\bar G}_{\varepsilon, \lambda}={\bar G}_c$, where ${\bar G}_{\varepsilon, \lambda}$ and ${\bar G}_c$ denotes the ${\mathcal T}_{\varepsilon, \lambda}$--  and ${\mathcal T}_c$-- closures of $G$, respectively.
\end{lemma}

\begin{lemma} Let $(E, {\mathcal P})$ be an $RLC$ module over $K$ with base $(\Omega, {\mathcal F}, P)$ such that $E$ has the countable concatenation property. Then the following are true:\\
\noindent (1). ${\bar G}_c={\bar G}_{\varepsilon, \lambda}$ has the countable concatenation property if so does $G$;\\
\noindent (2). If $G$ is $L^0$--convex, then ${\bar G}_{\varepsilon, \lambda}=[H_{cc}(G)]^-_{\varepsilon, \lambda}=[H_{cc}(G)]^-_c$ has the countable concatenation property.
\end{lemma}

\begin{proof} (1). $\bar{G}_c=\bar{G}_{\varepsilon,\lambda}$ is by Lemma 3.9, and thus we only need to prove that $\bar{G}_{\varepsilon,\lambda}$ has the countable concatenation property.

Let $\{x_n,n\in N\}$ be a given sequence in $\bar{G}_{\varepsilon,\lambda}$ and $\{A_n,n\in N\}$ a countable partition of $\Omega$ to $\mathcal{F}$, then by the countable concatenation property of $E$ there is $x^{\ast}\in E$ such that $x^{\ast}=\Sigma_{n=1}^{\infty}\tilde{I}_{A_n}x_n$. We claim that $x^{\ast}\in \bar{G}_{\varepsilon,\lambda}$, namely, $(x^{\ast}+N_{\theta}(\mathcal{Q},\varepsilon^{\ast},\lambda^{\ast}))\bigcap G\neq \emptyset$ for any given $\varepsilon^{\ast}>0,\lambda^{\ast}>0$ with $0<\lambda^{\ast}<1$ and any finite subset $\mathcal{Q}$ of $\mathcal{P}$, where $N_{\theta}(\mathcal{Q},\varepsilon^{\ast},\lambda^{\ast})=\{x\in E~|~P\{\omega\in\Omega~|~\|x\|_{\mathcal{Q}}(\omega)<\varepsilon^{\ast}\}>1-\lambda^{\ast}\}$.

In fact, it is clear that there exists $\bar{x}_n\in G$ for each $x_n\in\bar{G}_{\varepsilon,\lambda}$ such that $P\{\omega\in\Omega~|~\|x_n-\bar{x}_n\|_{\mathcal{Q}}(\omega)<\varepsilon^{\ast}\}>1-\frac{1}{2^{n+1}}\lambda^{\ast}$. By the countable concatenation property of $G$, there is $\bar{x}\in G$ such that $\bar{x}=\Sigma_{n=1}^{\infty}\tilde{I}_{A_n}\bar{x}_n$. Then $P\{\omega\in\Omega~|~\|x^{\ast}-\bar{x}\|_{\mathcal{Q}}(\omega)\geq\varepsilon^{\ast}\}=\Sigma_{n=1}^{\infty}P\{\omega\in A_n~|~\|x_n-\bar{x}_n\|_{\mathcal{Q}}(\omega)\geq\varepsilon^{\ast}\}\leq\Sigma_{n=1}^{\infty}P\{\omega\in \Omega~|~\|x_n-\bar{x}_n\|_{\mathcal{Q}}(\omega)\geq\varepsilon^{\ast}\}\leq\Sigma_{n=1}^{\infty}\frac{1}{2^{n+1}}\lambda^{\ast}=\frac{1}{2}\lambda^{\ast}$, namely $P\{\omega\in \Omega~|~\|x^{\ast}-\bar{x}\|_{\mathcal{Q}}(\omega)<\varepsilon^{\ast}\}\geq1-\frac{1}{2}\lambda^{\ast}>1-\lambda^{\ast}$, that is to say, $\bar{x}\in(x^{\ast}+N_{\theta}(\mathcal{Q},\varepsilon^{\ast},\lambda^{\ast}))\bigcap G$.

(2). By (1), $[H_{cc}(G)]^{-}_{c}=[H_{cc}(G)]^{-}_{\varepsilon,\lambda}$ has the countable concatenation property. Thus we only need to prove that $\bar{G}_{\varepsilon,\lambda}=[H_{cc}(G)]^{-}_{\varepsilon,\lambda}$. We can suppose , without loss of generality, that $\theta\in G$. Then for any $x=\Sigma_{n=1}^{\infty}\tilde{I}_{A_i}g_i\in H_{cc}(G)$, it is obvious that $\{\Sigma_{i=1}^{n}\tilde{I}_{A_i}g_i~|~n\in N\}$ is a $\mathcal{T}_{\varepsilon,\lambda}$--cauchy sequence in $G$ convergent to $x$ since $\{A_n,n\in N\}$ is a countable partition of $\Omega$ to $\mathcal{F}$, which means that $x\in \bar{G}_{\varepsilon,\lambda}$, namely $[H_{cc}(G)]^{-}_{\varepsilon,\lambda}\subset \bar{G}_{\varepsilon,\lambda}$, so $[H_{cc}(G)]^{-}_{\varepsilon,\lambda}=\bar{G}_{\varepsilon,\lambda}$.
\end{proof}

\begin{theorem} $(${\bf Random bipolar theorem}$)$ Let $\langle X, Y\rangle$ be a random duality pair over $K$ with base $(\Omega, {\mathcal F}, P)$ such that $X$ has the countable concatenation property. Then, for any subset $A$ of $X$, we have that $A^{00}=[H_{cc}(\Gamma (A))]^-_{\mathcal T}$ for each random compatible topology ${\mathcal T}$ with $Y$, where $\Gamma(A)$ denotes the $L^0$--balanced and $L^0$--convex hull of $A$ and $[H_{cc}(\Gamma (A))]^-_{\mathcal T}$ the ${\mathcal T}$--closure of $H_{cc}(\Gamma (A))$.
\end{theorem}

\begin{proof} Since $(X, {\mathcal T})^\ast=Y$, ${\mathcal T}\supset \sigma_c(X, Y)$. On the other hand, it is obvious that $A^{00}$ is an $L^0$--balanced, $L^0$--convex and $\sigma_c(X, Y)$--closed set with the countable concatenation property, so $A^{00}\supset [H_{cc}(\Gamma (A))]^-_{\sigma_c(X, Y)}\supset [H_{cc}(\Gamma (A))]^-_{\mathcal T}$. By (1) of Lemma 3.10 $[H_{cc}(\Gamma (A))]^-_{\mathcal T}$ has the countable concatenation property, if there is $x\in A^{00}\setminus [H_{cc}(\Gamma (A))]^-_{\mathcal T}$, then by Corollary 4.8 and Remark 4.9 of \cite{GZZ} there is $y\in (X, {\mathcal T})^\ast=Y$ such that $|\langle x, y\rangle|\nleqslant 1$ and $\bigvee\{|\langle a, y\rangle|\colon a\in A\}\leq 1$, which is impossible.
\end{proof}

\begin{remark} The classical bipolar theorem is an elegant result and hence frequently employed in the study of classical duality theory, cf. \cite{DXX}. However, the random bipolar theorem under the locally $L^0$--convex topology, namely Theorem 3.11 has the complicated form and also requires $X$ to have the countable concatenation property, so we do our best to avoid the use of it except in Subsection 3.3 where we are forced to use it to characterize a class of $L^0$--pre-barreled modules. In \cite{TXG-XXC} we proved a random bipolar theorem under the $(\varepsilon, \lambda)$--topology with the same shape as the classical bipolar theorem, but the countable concatenation property of $Y$ is required. To sum up, we are always forced to look for new methods in order to obtain some most refined results on random duality theory.
\end{remark}

It is time for us to speak of random compatible invariants. Corollary 4.8 of \cite{GZZ} shows that any closed $L^0$--convex sets with the countable concatenation property are random compatible invariants with respect to every random duality pair. Theorem 2.14 shows that the same is true for bounded sets in the sense of the locally $L^0$--convex topology.

\subsection{Random admissible topology}

\begin{definition} Let $\langle X, Y\rangle$ be a random duality pair over $K$ with base $(\Omega, {\mathcal F}, P)$ and ${\mathscr A}$ a family of $\sigma_c(Y, X)$--bounded sets of $Y$. For any $ A\in {\mathscr A}$, the $L^0$--seminorm $\|\cdot\|_A\colon X\to L^0_+({\mathcal F})$ is defined by $\|x\|_A=\bigvee\{|\langle x, a\rangle|\colon a\in A\}$ for all $x\in X$. Then the locally $L^0$--convex topology induced by the family $\{\|\cdot\|_A\colon A\in {\mathscr A}\}$ of $L^0$--seminorms, denoted by ${\mathcal T}_{\mathscr A}$, is called the topology of random uniform convergence on ${\mathscr A}$.
\end{definition}

\begin{definition} Let $\langle X, Y\rangle$ and ${\mathscr A}$ be the same as in Definition 3.13. ${\mathcal T}_{\mathscr A}$ is said to be random admissible if ${\mathcal T}_{\mathcal A}\supset \sigma_c(X, Y)$, in which case ${\mathscr A}$ is said to be random admissible. If $\mathcal{T}_{\mathscr{A}}$ is random compatible, namely $(X,\mathcal{T}_{\mathscr{A}})^{\ast}=Y$, then $\mathscr{A}$ is also said to be random compatible.
\end{definition}

As usual, let us first study $\mathcal{T}_{\mathscr{A}}$.

\begin{proposition} Let $\langle X,Y\rangle$ and $\mathscr{A}$ be the same as in Definition 3.13. Then the following are equivalent:

\noindent (1). $\mathcal{T}_{\mathscr{A}}$ is Hausdorff.

\noindent (2). $\bigcup\mathscr{A}:=\bigcup\{A:A\in\mathscr{A}\}$ is total, namely $\langle x,y\rangle=0$ for all $y\in \bigcup\mathscr{A}$ implies $x=\theta$.

\noindent (3). $Span\mathscr{A}$ (the submodule generated by $\bigcup\mathscr{A}$) is $\sigma_{\varepsilon,\lambda}(Y,X)$--dense in $Y$.

\noindent (4). $H_{cc}(Span\mathscr{A})$ is $\sigma_{c}(H_{cc}(Y),X)$--dense in $H_{cc}(Y)$.
\end{proposition}

\begin{proof} (1)$\Leftrightarrow$(2), (3)$\Rightarrow$(2) and (4)$\Rightarrow$(2) are all obvious.

(2)$\Rightarrow$(3). By Theorem 4.4 and Remark 4.9 of \cite{GZZ}, one can complete the proof by the same method as used in the classical case.

(3)$\Rightarrow$(4). By applying (2) of Lemma 3.10 to $Span\mathscr{A}$ and $(H_{cc}(Y),\sigma_{c}(H_{cc}(Y),X))$, we have the following relations:$$[H_{cc}(Span\mathscr{A})]^{-}_{\sigma_{c}(H_{cc}(X),Y)}$$$$=[H_{cc}(Span\mathscr{A})]^{-}_{\sigma_{\varepsilon,\lambda}(H_{cc}(X),Y)}$$
$$=[Span\mathscr{A}]^{-}_{\sigma_{\varepsilon,\lambda}(H_{cc}(X),Y)}~~~~~~~$$$$\supset H_{cc}(Y)~~~~~~~~~~~~~~~~~~~~~~~~$$
(by applying (3) to $H_{cc}(Y)$).
\end{proof}

Although random bipolar theorem does not necessarily hold for all random duality pairs, (2) of Lemma 3.16 below can complement this point.

\begin{lemma}\noindent Let $\langle X,Y\rangle$ be a random duality pair. Then we have:

\noindent (1). $A\subset Y$ is $\sigma_{c}(Y,X)$--bounded iff $A^{0}$ is a $\sigma_{c}(X,Y)$--$L^0$--barrel.

\noindent (2). For any $\sigma_{c}(Y,X)$--bounded set $A\subset Y$, $\|\cdot\|_B=\|\cdot\|_{B^{00}}$ (and hence $B^{00}$ is also $\sigma_{c}(Y,X)$--bounded), where $\|\cdot\|_B$ and $\|\cdot\|_{B^{00}}$ are defined as in Definition 3.13.
\end{lemma}

\begin{proof} (1) is clear.

(2). Since $B^{00}\supset B$, it is obvious that $\|\cdot\|_{B^{00}}\geq\|\cdot\|_B$. Conversely, if $\|x\|_B\leq 1$, then $x\in B^{0}$, and hence $\|x\|_{B^{00}}\leq 1$, which implies $\|\cdot\|_{B^{00}}\leq\|\cdot\|_B$.
\end{proof}

\begin{definition} Let $\langle X,Y\rangle$ be a random duality pair over $K$ with base $(\Omega,\mathcal{F},P)$ and $\mathscr{B}$ a family of $\sigma_{c}(Y,X)$--bounded sets of $Y$. $\mathscr{B}$ is saturated if the following are satisfied:

\noindent (a). If $A\subset B$ for some $B\in\mathscr{B}$ , then $A\in \mathscr{B}$;

\noindent (b). $A,B\in\mathscr{B}\Rightarrow A\bigcup B\in\mathscr{B}$;

\noindent (c). $B\in\mathscr{B}\Rightarrow B^{00}\in\mathscr{B}$;

\noindent (d). $\lambda B\in\mathscr{B}$ for all $\lambda\in L^0(\mathcal{F},K)$ and $B\in\mathscr{B}$.
\end{definition}

In the classical definition of a saturated family (which amouts to the case when $\mathcal{F}=\{\Omega,\emptyset\}$), the above (c) in Definition 3.17 is defined as ``$B\in\mathscr{B}\Rightarrow [\Gamma(B)]^{-}_{\sigma(Y,X)}\in\mathscr{B}$". But, generally, we only have the relation that $B^{00}\supset[\Gamma(B)]^{-}_{\sigma_c(Y,X)}$. Although the random bipolar theorem shows that $B^{00}=[H_{cc}(\Gamma(B))]^{-}_{\sigma_{c}(Y,X)}$ if $Y$ has the countable concatenation property, we would like to introduce the notion of a saturated family for an arbitrary random duality pair, so we choose Definition 3.17 to meet all our requirements.

Let $\langle X,Y\rangle$ be a random duality pair, in this paper we always denote by $\mathscr{B}(Y,X)$ the family of $\sigma_{c}(Y,X)$--bounded sets of $Y$ and $\beta(X,Y)=\mathcal{T}_{\mathscr{B}(Y,X)}$. By (2) of Lemma 3.16, $\mathscr{B}(Y,X)$ is saturated. It is also obvious that $\beta(X,Y)$ is the strongest random admissible topology.

For a family $\mathscr{A}$ of $\sigma_{c}(Y,X)$--bounded sets of $Y$, $\mathscr{A}^{s}$ denotes the saturated hull of $\mathscr{A}$, namely the smallest saturated family containing $\mathscr{A}$. It is easy to see that $\mathscr{A}^{s}=\{B\subset Y~|~$ there are $\lambda_1,\lambda_2,\cdots,\lambda_n\in L^{0}(\mathcal{F},K)$ and $A_1,A_2,\cdots,A_n\in\mathscr{A}$ such that $B\subset(\bigcup_{i=1}^{n}\lambda_i A_i)^{00}\}$, again by (2) of Lemma 3.16 one can easily see that $\mathcal{T}_{\mathscr{A}}=\mathcal{T}_{\mathscr{A}^s}$.

\begin{proposition} Let $\langle X,Y\rangle$ be a random duality pair over $K$ with base $(\Omega,\mathcal{F},P)$ and $\mathscr{A}$ and $\mathscr{B}$ two family of $\sigma_c(Y,X)$--bounded sets of $Y$ such that $\mathscr{B}$ is saturated. Then, $\mathcal{T}_{\mathscr{A}}\subset\mathcal{T}_{\mathscr{B}}$ iff $\mathscr{A}\subset\mathscr{B}$.
\end{proposition}

\begin{proof} If $\mathcal{T}_{\mathscr{A}}\subset\mathcal{T}_{\mathscr{B}}$, then for each $A\in\mathscr{A}$ there are $\xi\in L^{0}_{+}(\mathcal{F})$ and a finite subfamily $\{B_i~|~1\leq i\leq n\}$ of $\mathscr{B}$ such that $\|x\|_A\leq\xi(\bigvee\{\|x\|_{B_i}:1\leq i\leq n\})=\|x\|_B$ for all $x\in X$, where $B=\xi(\bigcup_{i=1}^{n}B_i)\in \mathscr{B}$. Thus $A\subset A^{00}\subset B^{00}\in\mathscr{B}$, which has showed that $A\in\mathscr{B}$. The converse is obvious.
\hfill $\square$
\end{proof}

\begin{corollary} Let $\langle X,Y\rangle$ be a random duality pair and $\mathscr{A}$ is a saturated family of $\sigma_{c}(Y,X)$--bounded sets of $Y$.  Then $\mathcal{T}_{\mathscr{A}}$ is random admissible iff $\bigcup\mathscr{A}=Y$.
\end{corollary}

\begin{proof} Let $Y_f$ denote the family of finite subsets of $Y$, then $\sigma_{c}(X,Y)=\mathcal{T}_{Y_f}$. So, $\mathcal{T}_{\mathscr{A}}$ is random admissible iff $\mathcal{T}_{Y_f}\subset\mathcal{T}_{\mathscr{A}}$ iff $Y_f\subset \mathscr{A}$ iff $\bigcup \mathscr{A}=Y$.
\end{proof}

\begin{theorem} Let $\langle X,Y\rangle$ be a random duality pair over $K$ with base $(\Omega,\mathcal{F},P)$. Then a locally $L^{0}$--convex topology $\mathcal{T}$ on $X$ is a topology of random uniform convergence iff $\mathcal{T}$ has a local base $\mathcal{B}$ at $\theta$ such that each $U\in\mathcal{B}$ is a $\sigma_c(X,Y)$--$L^{0}$--barrel ( where a $\sigma_c(X,Y)$--$L^{0}$--barrel means an $L^0$--balanced, $L^0$--absorbent, $L^0$--convex and $\sigma_c(X,Y)$--closed set of $X$).
\end{theorem}

\begin{proof}
If $\mathcal{T}=\mathcal{T}_{\mathscr{A}}$ for some family $\mathscr{A}$ of $\sigma_c(X,Y)-$bounded sets of $Y$, we can suppose that $\mathscr{A}$ is saturated, then $\mathcal{B}=\{A^{0}:A\in\mathscr{A}\}$ is a local base at $\theta$ of $\mathcal{T}$ such that each $A^{0}$ is a $\sigma_{c}(X,Y)$--$L^{0}$--barrel.

Conversely, let $\mathcal{T}$ have a local base $\mathcal{B}$ at $\theta$ such that each $U\in \mathcal{B}$ is a $\sigma_c(X,Y)$--$L^0$--barrel. Let $\mathscr{A}=\{U^{0}:U\in\mathscr{B}\}$, then each $U^{0}$ is $\sigma_{c}(Y,X)$--bounded since $(U^{0})^{0}\supset U$ is $L^0-$absorbent. Further, we show that $\mathcal{T}=\mathcal{T}_{\mathscr{A}}$ as follows. Since $\mathcal{T}_{\mathscr{A}}$ is induced by $\{\|\cdot\|_{A}:A\in\mathscr{A}\}$ and , for each $U^{0}\in\mathscr{A}$, $\{x\in X~|~\|x\|_{U^{0}}\leq 1\}=U^{00}\supset U$, which shows that $\|\cdot\|_{U^0}$ is $\mathcal{T}$--continuous, namely $\mathcal{T}_{\mathscr{A}}\subset \mathcal{T}$. On the other hand, for each $U\in\mathcal{B}$, $U\subset U^{00}=(U^{0})^{0}$, namely each element $(U^{0})^{0}$ of a local base at $\theta$ of $\mathcal{T}_{\mathscr{A}}$ is a $\mathcal{T}$--neighborhood of $\theta$, that is to say, $\mathcal{T}_{\mathscr{A}}\subset\mathcal{T}$.
\end{proof}

In the classical case, by the classical bipolar theorem it can be easily established that $\{A^{0}:A\in\mathscr{B}(Y,X)\}$ as the local base at $\theta$ of $\beta(X,Y)$ is exactly the family of $\sigma(X,Y)$--barrels. However, in the random setting, we do not know if $\{A^{0}:A\in\mathscr{B}(Y,X)\}$ as the local base at $\theta$ of $\beta(X,Y)$ is still the family of $\sigma_{c}(X,Y)$--$L^0$--barrels, we only know that for each $\sigma_{c}(X,Y)$--$L^0$--barrel $U$ there is $A(=U^0)\in\mathscr{B}(Y,X)$ such that $U\subset A^{0}$. So, we remind the reader of the following useful result:

\begin{theorem} Let $\langle X,Y\rangle$ be a random duality pair such that $X$ has the countable concatenation property. Then the family of  $\sigma_c(X,Y)$--$L^0$--barrels with the countable concatenation property forms a local base at $\theta$ of $\beta(X,Y)$.
\end{theorem}

\begin{proof} By the countable concatenation property of $X$, it is easy to see that $A^0$ has the countable concatenation property for each $A\in\mathscr{B}(Y,X)$, and hence $A^0$ is a $\sigma_{c}(X,Y)$--$L^0$--barrel with the countable concatenation property. On the other hand, for each $\sigma_{c}(X,Y)$--$L^0$--barrel $U$ with the countable concatenation property, then by Theorem 3.11 we have that $U=U^{00}=(U^{0})^{0}$. Since $U^{0}\in\mathscr{B}(Y,X)$, $U\in\{A^{0}:A\in\mathscr{B}(Y,X)\}$. To sum up, the family of $\sigma_{c}(X,Y)$--$L^{0}$--barrels with the countable concatenation property is exactly the local base $\{A^{0}:A\in\mathscr{B}(Y,X)\}$.
\end{proof}

Theorem 3.22 below shows that the study of random admissible topology is of universal interest in the theory of $RCL$ modules.

\begin{theorem} Let $(X,\mathcal{P})$ be an $RLC$ module over $K$ with base $(\Omega,\mathcal{F},P)$ and $\mathcal{E}$ the family of all the subsets $E$ of $X_{c}^{\ast}$ such that $E$ is equicontinuous from $(X,\mathcal{T}_c)$ to $L^{0}(\mathcal{F},K)$ endowed with the locally $L^{0}$--convex topology induced by $|\cdot|$. Then $\mathcal{T}_{c}=\mathcal{T}_{\mathcal{E}}$, where we consider the natural pairing $\langle X,X_{c}^{\ast}\rangle$, then $\mathcal{T}_{\mathcal{E}}$ is, clearly, a random admissible topology.
\end{theorem}

\begin{proof} It is clear that $E\in\mathcal{E}$ iff there are $\xi\in L^{0}_{+}(\mathcal{F})$ and a finite subset $\mathcal{Q}$ of $\mathcal{P}$ such that $\|x\|_{E}:=\bigvee\{|f(x)|:f\in E\}\leq \xi(\bigvee\{\|x\|:\|\cdot\|\in\mathcal{Q}\})$ for all $x\in X$, so $\mathcal{T}_{\mathcal{E}}\subset\mathcal{T}_{c}$.

Conversely, for each $\|\cdot\|\in\mathcal{P}$, let $E=\{f\in X^{\ast}_{c}~|~|f(x)|\leq \|x\|$ for all $x\in X\}$, then from the random Hahn-Banach theorem of \cite{TXG-JFA} one can easily see that $\|\cdot\|=\|\cdot\|_E$, so $\mathcal{T}_{c}\subset\mathcal{T}_{\mathcal{E}}$.
\end{proof}

\begin{corollary} Let $\langle X,Y\rangle$ be a random duality pair over $K$ with base $(\Omega,\mathcal{F},P)$. Then every random compatible topology  $\mathcal{T}$ on $X$ is random admissible.
\end{corollary}

\begin{proof} By Definition 3.2, there is a family $\mathcal{P}$ of $L^0$--seminorms on $X$ such that $(X,\mathcal{P})$ becomes an $RLC$ module over $K$ with base $(\Omega,\mathcal{F},P)$ and $(X,\mathcal{P})^{\ast}_{c}=Y$, $\mathcal{T}$ is just induced by $\mathcal{P}$, at which time $\langle X,Y\rangle$ is exactly $\langle X,X^{\ast}_{c}\rangle$ and $\mathcal{T}=\mathcal{T}_{\mathcal{E}}$ by Theorem 3.22.
\end{proof}

The proof of Theorem 3.24 below (namely the resonance theorem) is omitted since it is the same as that of the classical case.
For the notions of an $L^0$--barreled module and $L^0$--pre--barreled module, see the beginning part of Subsection 3.3 of this paper.

\begin{theorem} Let $(E,\mathcal{P})$ be an $RLC$ module over $K$ with base $(\Omega,\mathcal{F},P)$ and $H\subset E^{\ast}_{c}$. Then we have the following:

\noindent $(1)$. If $(E,\mathcal{T}_c)$ is $L^0$--barreled, then $H$ is equicontunuous from $(E,\mathcal{T}_c)$ to $(L^{0}(\mathcal{F},K),\mathcal{T}_c)$ iff $H$ is $\sigma_c(E^{\ast}_{c},E)$--bounded.

\noindent $(2)$. If $(E,\mathcal{T}_c)$ is $L^0$--pre-barreled and $E$ has the countable concatenation property, then $H$ is equicontinuous from $(E,\mathcal{T}_c)$ to $(L^{0}(\mathcal{F},K),\mathcal{T}_c)$ iff $H$ is $\sigma_c(E^{\ast}_{c},E)$--bounded.
\end{theorem}

In the classical case, for a locally convex space $(E,\mathcal{T})$, a subset $H\subset E^{\ast}$ is equicontinuous, then it must be $\sigma(E^{\ast},E)$--relatively compact. However the classical Banach-Alaoglu theorem universally fails to hold in the case of $RN$ modules under the $(\varepsilon,\lambda)$--topology (cf. \cite{TXG-Alao}), the same, of course, occurs for the locally $L^{0}$--convex topology, so we can not generalize the construction of the classical Mackey topology to the random setting.

\subsection{A characterization for a random locally convex module to be $L^0$--pre-barreled}

Let us first recall the notion of an $L^0$--barreled module from \cite{FKV}. Let $(E, {\mathcal T})$ be a locally $L^0$--convex $L^0({\mathcal F}, K)$--module. An $L^0$--balanced, $L^0$--absorbent, $L^0$--convex and closed subset of $E$ is called an $L^0$--barrel. $(E, {\mathcal T})$ is called an $L^0$--barreled module if every $L^0$--barrel is a neighborhood of $\theta$, whereas it seems to us that the following notion of an $L^0$--pre-barreled module is more suitable for financial applications.

\begin{definition} Let $(E, {\mathcal T})$ be a locally $L^0$--convex $L^0({\mathcal F}, K)$--module. $(E,$ ${\mathcal T})$ is called an $L^0$--pre-barreled module if every $L^0$--barrel with the countable concatenation property is a neighborhood of $\theta$.
\end{definition}

The main result of this subsection is Theorem 3.26 below.

\begin{theorem}  Let $(E,\mathcal{P})$ be an $RLC$ module over $K$ with base $(\Omega,\mathcal{F},P)$ such that $E$ has the countable concatenation property. Then $(E,\mathcal{T}_c)$ is $L^0$--pre-barreled iff $\mathcal{T}_c=\beta(E,E^{\ast}_c)$.
\end{theorem}

\begin{proof} $\mathcal{T}_c$ has a local base at $\theta$ consisting of $\{N_{\theta}(\mathcal{Q},\varepsilon)~|~\mathcal{Q}\in\mathcal{P}_f,\varepsilon\in L^{0}_{++}(\mathcal{F})\}$, where $N_{\theta}(\mathcal{Q},\varepsilon)=\{x\in E~|~\|x\|_{\mathcal{Q}}\leq\varepsilon\}$. It is obvious that every $N_{\theta}(\mathcal{Q},\varepsilon)$ is an $\mathcal{T}_c$--$L^0$--barrel with the countable concatenation property. Thus, if $\mathcal{T}_c$ is $L^0$--pre-barreled, then the family of $\mathcal{T}_c$--$L^0$--barrels with the countable concatenation property forms a local base at $\theta$ of $\mathcal{T}_c$. By Corollary 4.8 of \cite{GZZ}, $L^0$--barrels with the countable concatenation property are random compatible invariants, namely the family of $\mathcal{T}_c$--$L^0$--barrels with the countable concatenation property coincides with the family of $\sigma_c(E,E_c^{\ast})$--$L^0$--barrels with the countable concatenation property, so $\mathcal{T}_c=\beta(E,E_c^{\ast})$ by Theorem 3.21 if $\mathcal{T}_c$ is $L^0$--pre-barreled.

Conversely, if $\mathcal{T}_c=\beta(E,E_c^{\ast})$, then ,since every $\mathcal{T}_c$--$L^0$--barrel with the countable concatenation property is $\sigma_c(E,E_c^{\ast})$--$L^0$--barrel by Corollary 4.8 of \cite{GZZ}, and hence a $\beta(E,E^{\ast}_c)$--neighborhood of $\theta$ by Theorem 3.21, namely a $\mathcal{T}_c$--neighborhood of $\theta$, that is to say, $(E,\mathcal{T}_c)$ is $L^0$--pre-barreled.
\end{proof}

\begin{corollary} Let $(E,\|\cdot\|)$ be a $\mathcal{T}_c$--complete $RN$ module over $K$ with $(\Omega,\mathcal{F},P)$ such that $E$ has the countable concatenation property. Then $(E,\mathcal{T}_c)$ is $L^0$--pre-barreled.
\end{corollary}

\begin{proof}  We only need to verify that $\mathcal{T}_c=\beta(E,E_c^{\ast})$.

First, the locally $L^0$--convex topology $\mathcal{T}_c$ induced by $\|\cdot\|$ is a random compatible topology with respect to the natural random duality pair $\langle E,E_c^{\ast}\rangle$, so $\mathcal{T}_c\subset\beta(E,E_c^{\ast})$ by Corollary 3.23.

Conversely, $\beta(E,E_c^{\ast})$ is induced by $\{\|\cdot\|_A:A\in \mathscr{B}(E_c^{\ast},E)\}$, please recall that $\|\cdot\|_A:E\rightarrow L^{0}_{+}(\mathcal{F})$ is given by $\|x\|_A=\bigvee\{|f(x)|:f\in A\}$ for all $x\in E$ and $A\in \mathscr{B}(E_c^{\ast},E)$. Thus we only need to prove that each $\|\cdot\|_A$ is continuous from $(E,\mathcal{T}_c)$ to $(L^0(\mathcal{F},K),\mathcal{T}_c)$. $A\in \mathscr{B}(E_c^{\ast},E)$ means that $\{f(x):f\in A\}$ is $\mathcal{T}_c$--bounded in $(L^0(\mathcal{F},K),\mathcal{T}_c)$ for each $x\in E$, then, by Theorem 2.13 $A$ is $\mathcal{T}_c$--bounded in $E_c^{\ast}$, namely a.s. bounded, and hence there is $\xi_A\in L^0_+(\mathcal{F})$ such that $\|f\|\leq\xi_A$ for all $f\in A$. This shows that $\|x\|_A=\bigvee\{|f(x)|:f\in A\}\leq\bigvee\{\|f\|\cdot\|x\|:f\in A\}\leq(\bigvee\{\|f\|:f\in A\})\|x\|\leq\xi_A\|x\|$ for all $x\in E$, namely $\beta(E,E_c^{\ast})\subset \mathcal{T}_c$. \hfill $\square$
\end{proof}

\begin{corollary} For each $p\in[1,+\infty]$, $(L^{p}_{\mathcal{F}}(\mathcal{E}),\mathcal{T}_c)$ is $L^0$--pre-barreled.
\end{corollary}

\begin{proof} Since $L^{p}_{\mathcal{F}}(\mathcal{E})$ is $\mathcal{T}_c$--complete and has the countable concatenation property, then it immediately follows from Corollary 3.27.
\end{proof}

\section{Continuity and subdifferentiability theorems in $L^0$--pre-barreled modules}

We will state the results in this section under the framework of locally $L^0$--convex modules since the proofs of these results do not necessarily depend on the family of $L^0$--seminorms. Continuity and subdifferentiability theorems in $L^0$--barreled modules were already proved in \cite{FKV}. As shown in \cite{FKV}, the proofs in the random setting are very similar to those in the corresponding classical cases. Thus this section  is focused on some discussions on the relation between the topological structure and stratification structure of a locally $L^0$--convex module.

In the section, a locally $L^0$--convex module (a topological $L^0$--module) means a locally $L^0$--convex $L^0({\mathcal{F},R})$--module (resp., a topological $L^0({\mathcal{F},R})$--module).

To state our main results, let us first recall the following: let $E$ be an $L^0({\mathcal F})$--module and $f$ a function from $E$ to ${\bar L}^0({\mathcal F})$. The effective domain of $f$ is denoted by $dom(f):=\{x\in E~|~f(x)<+\infty~\hbox{on}~\Omega\}$ and the epigraph of $f$ by $epi(f):=\{(x,r)\in E\times L^0({\mathcal F})~|~f(x)\leq r\}$. $f$ is proper if $dom(f)\neq\emptyset$ and $f(x)>-\infty~\hbox{on}~\Omega$. $f$ is $L^0$--convex if $f(\xi x+(1-\xi)y)\leq \xi f(x)+(1-\xi)f(y)$ for all $x,~y\in E$ and $\xi\in L^0_+({\mathcal F})$ with $0\leq \xi\leq 1$, where the following convention is adopted: $0\cdot(\pm\infty)=0$ and $+\infty\pm(\pm\infty)=+\infty$. $f\colon E\to {\bar L}^0({\mathcal F})$ is said to be local ( or, to have the local property ) if ${\tilde I}_Af(x)={\tilde I}_Af({\tilde I}_Ax)$ for all $x\in E$ and $A\in {\mathcal F}$. In \cite{FKV-appro}, it is proved that an $L^0$-convex function is local. In this paper, we adopt the weakest definition of lower semicontinuity, namely, let $(E,\mathcal{T})$ be a locally $L^0$--convex module and $f:E\rightarrow\bar{L}^0(\mathcal{F})$ a proper $L^0$--convex function, we say that $f$ is lower semicontinuous if $\{x\in E~|~f(x)\leq r\}$ is ${\mathcal T}$-closed for any $r\in L^0({\mathcal F})$.

\begin{definition}(See \cite{FKV}). Let $(E,\mathcal{T})$ be a locally $L^0$--convex module and $f:E\rightarrow\bar{L}^0(\mathcal{F})$ a proper $L^0$--convex function. $f$ is subdifferentiable at $x\in dom(f)$ if there is $u\in (E,\mathcal{T})^\ast$ such that $u(y-x)\leq f(y)-f(x)$ for all $y\in E$, at which time $u$ is called a subgradient of $f$ at $x$. The set of subgradients of $f$ at $x$ is denoted by $\partial f(x)$.
\end{definition}

We can now state our main results as follows:

\begin{theorem} Let $(E,\mathcal{T})$ be a real $L^0$--pre-barreled module such that $E$ has the countable concatenation property. Then a proper lower semicontinuous $L^0$--convex function $f:E\rightarrow \bar{L}^{0}(\mathcal{F})$ is continuous on $Int(dom(f)):=$ the interior of $dom(f)$, namely $f$ is continuous from $(Int(dom(f)),\mathcal{T})$ to $(L^0(\mathcal{F}),\mathcal{T}_c)$.
\end{theorem}

\begin{theorem} Let $(E,\mathcal{T})$ be a real $L^0$--pre-barreled module such that $E$ has the countable concatenation property. Then, for a proper lower semicontinuous $L^0$--convex function $f:E\rightarrow\bar{L}^{0}(\mathcal{F})$, $\partial f(x)\neq\emptyset$ for all $x\in Int(dom(f))$.
\end{theorem}

To prove Theorem 4.2, we needs the following known lemmas:

\begin{lemma} ($See$ \cite{FKV}). Let $E$ be a topological $L^0$--module. If in some neighborhood of an element $x_{0}\in E$ a proper $L^0$--convex function $f:E\rightarrow\bar{L}^{0}(\mathcal{F})$ is bounded above by some $\xi_0\in L^{0}(\mathcal{F})$, then $f$ is continuous at $x_{0}$.
\end{lemma}

\begin{lemma} ($See$ \cite{FKV}). Let $E$ be a topological $L^0$--module and $f:E\rightarrow\bar{L}^{0}(\mathcal{F})$ a proper $L^0$--convex function. Then the following statements are equivalent:

\noindent $(1)$. There is a nonempty open set $O\subset E$ on which $f$ is bounded above by some $\xi_0\in L^{0}(\mathcal{F})$.

\noindent $(2)$. $f$ is continuous on $Int(dom(f))$ and $Int(dom(f))\neq\emptyset$.
\end{lemma}

\begin{lemma} ($See$ \cite{FKV}). Let $E$ be a topological $L^0$--module and $x\in E$. Then every proper $L^0$--convex function $f:Span_{L^0}(x)\rightarrow \bar{L}^{0}(\mathcal{F})$ is continuous on $Int(dom(f))$, where $Span_{L^0}(x)$ is the $L^0$--module spanned by $x$ and endowed with the relative topology.
\end{lemma}

We can now prove Theorem 4.2.

\noindent{\bf Proof of Theorem 4.2.}
Assume that there is $x_{0}\in Int(dom(f))$. By translation, we may assume $x_{0}=0$ and further take $Y_0\in L^{0}(\mathcal{F})$ such that $f(0)<Y_0$ on $\Omega$. Since $f$ is lower semicontinuous, the set $C:=\{x\in E~|~f(x)\leq Y_0\}$ is closed. Further, for all $x\in E$, the net $\{\frac{x}{Y}:Y\in L^{0}_{++}(\mathcal{F})\}$ converge to $\theta$. By Lemma 4.6, the restriction of $f$ to $Span_{L^0}(x)$ is continuous at $\theta$, hence $f(\frac{x}{Y})<Y_0$ on $\Omega$ for large $Y$, which means that $C$ is $L^0$--absorbent. Hence $C\bigcap(-C)$ is an $L^{0}$--barrel. Since $E$ has the countable concatenation property and $f$ has the local property, it is easy to observe that $C\bigcap(-C)$ is an $L^0$--barrel with the countable concatenation property and in turn a neighborhood of $\theta\in E$, so $f$ is continuous on $Int(dom(f))$ by Lemma 4.5.   \hfill $\square$

To prove Theorem 4.3, we need the following three lemmas as well as Theorem 4.10 below.

Lemma 4.7 below is a slight generalization of Lemma 3.17 of \cite{TXG-JFA}, whereas their proofs are the same, so the proof of Lemma 4.7 is omitted.

\begin{lemma} Let $(E,\mathcal{T})$ be a locally $L^0$--convex $L^0(\mathcal{F},K)$--module and $A\in \cal{F}$ with $P(A)>0$. If $G$ and $M$ are an open set and a closed set of $E$, respectively, such that $\tilde{I}_A G+\tilde{I}_{A^{c}} G\subset G$ and $\tilde{I}_A M+\tilde{I}_{A^{c}} M\subset M$, then $\tilde{I}_A G$ is relatively open in $\tilde{I}_A E$ and $\tilde{I}_A M$ is relatively closed in $\tilde{I}_A E$.
\end{lemma}

From Lemma 4.7 one can see that $(\tilde{I}_{A}E,\mathcal{T}|_{\tilde{I}_{A}E})$ is still a locally $L^0$--convex $L^0(\mathcal{F}_A,K)$--module, where $\mathcal{F}_A=A\bigcap\mathcal{F}:=\{A\bigcap B~|~B\in\mathcal{F}\}$ is the $\sigma$--algebra of $(A,\mathcal{F}_A,P(\cdot|A))$.

\begin{lemma} Let $(E,\mathcal{T})$ be a locally $L^0$--convex module, $A\in\mathcal{F}$ with $P(A)>0$ and $f:E\rightarrow \bar{L}^0(\mathcal{F})$ a proper $L^0$--convex function. If, we define $f_A:\tilde{I}_{A}E\rightarrow \tilde{I}_{A}\bar{L}^0(\mathcal{F})$ by $f_A(\tilde{I}_{A}x)=\tilde{I}_{A}f(\tilde{I}_{A}x)$ for all $x\in E$, then we have:

\noindent $(1)$. For all $x\in dom(f)$, $\tilde{I}_A(x,f(x))\in \partial(epi(f_A))$, where $\partial(epi(f_A))$ denotes the boundary of $epi(f_A)$ in $(\tilde{I}_{A}E,\mathcal{T}|_{\tilde{I}_{A}E})\times(\tilde{I}_{A}L^0(\mathcal{F}),\mathcal{T}_c|$ $_{\tilde{I}_A L^0{(\mathcal{F})}})$. Here, let us recall that $\mathcal{T}_c$ is the locally $L^0$--convex topology on $L^0(\mathcal{F})$ induced by $|\cdot|$.

\noindent $(2)$. For all $x\in dom(f)$, $\tilde{I}_{A}(x,f(x))\not\in \tilde{I}_{A}(Int(epi(f)))$.
\end{lemma}

\begin{proof} (1). It is easy to see that $(x,f(x))\in\partial(epi(f))$ for all $x\in dom(f)$. By the local property of $f$, $\tilde{I}_{A}(x,f(x))=(\tilde{I}_{A}x, \tilde{I}_{A}f(x))=(\tilde{I}_{A}x,f_A(\tilde{I}_{A}x))$ for all $x\in E$. So, if we consider the corresponding problem in $(\tilde{I}_{A}E,\mathcal{T}|_{\tilde{I}_{A}E})$, then we have that $\tilde{I}_{A}(x,f(x))=(\tilde{I}_{A}x,f_A(\tilde{I}_{A}x))\in\partial(epi(f_A))$ for all $x\in dom(f)$.

(2). By the above (1), it is , of course, that $\tilde{I}_{A}(x,f(x))\not\in Int_A(epi(f_A))$, where $Int_A(epi(f_A))$ denotes the interior of $epi(f_A)$ in $(\tilde{I}_{A}E,\mathcal{T}|_{\tilde{I}_{A}E})\times(\tilde{I}_{A}L^0(\mathcal{F}),$ $\mathcal{T}_c|_{\tilde{I}_A  L^0{(\mathcal{F})}})$. It is obvious that $epi(f_A)=\tilde{I}_A(epi(f))$. By Lemma 4.7, $\tilde{I}_{A}(Int(epi(f)))$ is an open set in $(\tilde{I}_{A}E,\mathcal{T}|_{\tilde{I}_{A}E})\times(\tilde{I}_{A}L^0(\mathcal{F}),$ $\mathcal{T}_c|_{\tilde{I}_A L^0{(\mathcal{F})}})$, so $Int_A(epi(f_A))=Int_A(\tilde{I}_A epi(f))\supset \tilde{I}_A(Int(epi(f)))$, which implies that $\tilde{I}_A(x,$ $f(x))\not\in \tilde{I}_A(Int$ $(epi(f)))$ for all $x\in dom(f)$. \hfill $\square$
\end{proof}

Proof of Lemma 4.9 below is the same as that of Lemma 3.14 of \cite{FKV}, so is omitted.

\begin{lemma} Let $(E,\mathcal{T})$ be a locally $L^{0}$--convex module and $f:E\rightarrow\bar{L}^{0}(\mathcal{F})$ a proper lower semicontinuous $L^{0}$--convex function. Then $Int(epi(f))\neq\emptyset$ implies $Int(dom(f))\neq\emptyset$. Furthermore, if, in addition, $(E,\mathcal{T})$ is $L^0$--pre-barreled such that $E$ has the countable concatenation property, then $Int$ $(dom(f))$ $\neq\emptyset$ iff $Int(epi(f))\neq\emptyset$.
\end{lemma}

\begin{theorem}($See$ \cite{FKV}). Let $(E, {\mathcal T})$ be a Hausdorff locally $L^0$--convex $L^0({\mathcal F}, K)$--module,
$M$ and $G$ two nonempty $L^0$--convex sets of $E$ with $G$ open. If ${\tilde I}_AM\cap {\tilde I}_AG=\emptyset$ for all $A\in {\mathcal F}$ with $P(A)>0$, then there is $f\in (E,\mathcal{T})^\ast$ such that:$$(Ref)(y)<(Ref)(z)~\hbox{on $\Omega$ for all $y\in G$ and $z\in M$}.$$
\end{theorem}

We can now prove Theorem 4.3.

\noindent{\bf Proof of Theorem 4.3.}
Let $x_0\in Int(dom(f))$. We separate $(x_0,f(x_0))$ from $Int(epi(f))$ by means of Theorem 4.10. By Lemma 4.9, $Int(epi(f))$ is not-empty, $(x_0,f(x_0))\in \partial(epi(f))$ and $$\tilde{I}_{A}(x_0,f(x_0))\bigcap\tilde{I}_{A}Int(epi(f))=\emptyset$$ for all $A\in\mathcal{F}$ and $P(A)>0$. Hence, there are $g_1\in (E,\mathcal{T})^\ast$ and $g_2\in (L^{0}(\mathcal{F}), {\mathcal T}_c)^{\ast}$ (in fact, $(L^{0}(\mathcal{F}), {\mathcal T}_c)^{\ast}=L^{0}(\mathcal{F})$) such that $$g_1(x)+g_{2}(y)<g_1(x_0)+g_2(f(x_0))$$ on $\Omega$ for all $(x,y)\in Int(epi(f))$. By the fact that $g_2(y)=yg_2(1)$ we derive that $g_2(1)<0$ on $\Omega$. We will show that $\frac{-g_1}{g_2(1)}\in \partial f(x_0)$. To this end, let $x\in E$, $A=[f(x)=+\infty]$ and $\tilde{x}=I_{A}x_0+I_{A^c}x$. Then, $\tilde{x}\in dom(f)$ and in turn $(\tilde{x},f(\tilde{x}))\in \partial(epi(f))$. Thus, there is a net $(x_{\alpha},y_{\alpha})\subset Int(epi(f))$ which converges to $(\tilde{x},f(\tilde{x}))$ and for which \\
$(1) ~~~~~~~~~~~~~~~~~~~~g_1(x_{\alpha})+y_{\alpha}g_2(1)<g_1(x_0)+g_2(f(x_0))$\\
on $\Omega$ for all $\alpha$. Since $g_1\in (E,\mathcal{T})^\ast$ we may pass to limits in (1) yielding $$\frac{-g_1(\tilde{x}-x_0)}{g_2(1)}\leq f(\tilde{x})-f(x_0).$$ Finally, from the local property of $f$ and $g_1$ we derive $$\frac{-g_1(\tilde{x}-x_0)}{g_2(1)}\leq f(x)-f(x_0)$$ and since $x\in E$ is arbitrary we conclude that $\frac{-g_1}{g_2(1)}$ indeed is subgradient of $f$ at $x_0$.  \hfill $\square$

If the hypothesis ``that $(E,\mathcal{T})$ is $L^0$--pre-barreled module such that $E$ has the countable concatenation property" in Theorems 4.2 and 4.3 is replaced by the one ``that $(E,\mathcal{T})$ is $L^0$--barreled", then Theorems 4.2 and 4.3 change to Proposition 3.5 and Theorem 3.7 of \cite{FKV}, respectively, and the proofs of Theorems 4.2 and 4.3 are also very similar to that of Proposition 3.5 and Theorem 3.7 of \cite{FKV}, what is different is that Theorems 4.2 and 4.3 are based on the notion of $L^0$--pre-barreled modules, at least Theorems 4.2 and 4.3 make sense for $L^0$--pre-barreled $RLC$ modules with the countable concatenation property since we have characterized such $RLC$ modules in Theorem 3.26, furthermore, Corollary 3.28 shows that Corollaries 4.11 and 4.12 below of Theorems 4.2 and 4.3, respectively, are enough for financial applications, for whose proofs one only notices that $(E,\mathcal{T}_c)$ is a Hausdorff locally $L^0$--convex module with $\mathcal{T}_c$ induced by $\mathcal{P}$ for a random locally convex module $(E,\mathcal{P})$ over $R$ with base $(\Omega,\mathcal{F},P)$.

\begin{corollary} Let $(E,\mathcal{P})$ be an $L^0$--pre-barreled random locally convex module over $R$ with base $(\Omega,\mathcal{F},P)$ such that $E$ has the countable concatenation property. Then a proper $\mathcal{T}_c$--lower semicontinuous $L^0$--convex function $f:E\rightarrow \bar{L}^{0}(\mathcal{F})$ is continuous on $Int(dom(f)):=$ the $\mathcal{T}_c$--interior of $dom(f)$, namely $f$ is continuous from $(Int(dom(f)),\mathcal{T}_c)$ to $(L^0(\mathcal{F}),\mathcal{T}_c)$.
\end{corollary}

\begin{corollary} Let $(E,\mathcal{P})$ be an $L^0$--pre-barreled random locally convex module over $R$ with base $(\Omega,\mathcal{F},P)$ such that $E$ has the countable concatenation property. Then, for a proper $\mathcal{T}_c$--lower semicontinuous $L^0$--convex function $f:E\rightarrow\bar{L}^{0}(\mathcal{F})$, $\partial f(x)\neq\emptyset$ for all $x\in Int(dom(f))$, where $Int(dom(f))$ is the same as in Corollary 4.11.
\end{corollary}

\sec{Acknowledgements} {This work was supported by National Natural Science Foundation of China (Grant No. 11171015 and 11301568).}


\end{document}